\documentclass[11pt]{article}

\usepackage{amssymb,bbm,comment}
\usepackage{amsthm, amsmath}

\renewcommand{\cal}{\mathcal}   

\renewcommand{\mathbb}{\mathbbm}                   
\renewcommand{\epsilon}{\varepsilon}
\renewcommand{\phi}{\varphi}
\renewcommand{\theta}{\vartheta}
\renewcommand{\le}{\leqslant}
\renewcommand{\ge}{\geqslant}

\newcommand{\law}{{\cal L}}
\newcommand{\abs}[1]{\left\lvert #1 \right\rvert}   
\newcommand{\babs}[1]{\Bigl\lvert #1 \Bigr\rvert}   
\newcommand{\scapro}[2]{\langle #1,#2\rangle}       
\newcommand{\norm}[1]{\left\lVert #1 \right\rVert}  

\DeclareMathOperator{\1}{\mathbbm 1}                 
\DeclareMathOperator{\E}{{\mathbb E}}               
\DeclareMathOperator{\Var}{Var}                     
\DeclareMathOperator{\Cov}{Cov}                     
\DeclareMathOperator{\R}{{\mathbb R}}               
\DeclareMathOperator{\Rp}{{\mathbb R}_+}            
\DeclareMathOperator{\C}{{\mathbb C}}               
\DeclareMathOperator{\N}{{\mathbb N}}               
\DeclareMathOperator{\Id}{Id}                       
\DeclareMathOperator{\F}{{\mathcal F}} \DeclareMathOperator{\calP}{\mathcal{P}}
\DeclareMathOperator{\Borel}{{\mathcal B}}
\renewcommand{\Re}{\text{Re}}               
\DeclareMathOperator{\TV}{TV}

\numberwithin{equation}{section} 

{\theoremstyle{plain}
\newtheorem{definition}{Definition}[section]
\newtheorem{lemma}[definition]{Lemma}
\newtheorem{theorem}[definition]{Theorem}
\newtheorem{corollary}[definition]{Corollary}
\newtheorem{proposition}[definition]{Proposition}
\newtheorem{remark}[definition]{Remark}
\newtheorem{examples}[definition]{Examples}
\newtheorem{assumption}[definition]{Assumption}
}

\begin{document}

\title{Delay differential equations driven by L\'evy processes:
stationarity and Feller properties}



\author{Markus Rei{\ss}\\
   Weierstra{\ss} Institute for\\ 
   Applied Analysis and Stochastics\\
  Mohrenstrasse 39\\10117 Berlin, Germany\\
\and Markus Riedle \\
    Institute of Mathematics\\
    Humboldt-University of Berlin\\
    Unter den Linden 6\\
    10099 Berlin, Germany\\
\and Onno van Gaans\thanks{ O.\ van\ Gaans acknowledges the financial support provided
through the European Community's Human Potential Programme under contracts
HPRN-CT-2000-00100, DYNSTOCH and
HPRN-CT-2002-00281.}\\
     Department for Mathematics\\ and Computer Science\\
     Friedrich Schiller University Jena\\
     07740 Jena, Germany\\
 }

\date{}

\maketitle

\begin{abstract}
We consider a stochastic delay differential equation driven by a general L\'evy process.
Both, the drift and the noise term may depend on the past, but only the drift term is
assumed to be linear. We show that the segment process is eventually Feller, but in
general not eventually strong Feller on the Skorokhod space. The existence of an
invariant measure is shown by proving tightness of the segments using semimartingale
characteristics and the Krylov-Bogoliubov method. A counterexample shows that the
stationary solution in completely general situations may not be unique, but in more
specific cases uniqueness is established.
\end{abstract}



\section{Introduction}

Stochastic delay differential equations, also known as stochastic functional
differential
equations, are a natural generalisation of stochastic ordinary differential
equations by
allowing the coefficients to depend on values in the past.  When only the drift
coefficient depends on the past, main stochastic tools and results for
stochastic
ordinary differential equations can be applied, for example by removing the
drift via a
change of measure. If the stochastic perturbation depends on the past, however,
surprising new phenomena emerge, see Mohammed and Scheutzow \cite{MohSch97} for
a
discussion on flow and stability properties.

The main purpose of the present work is to investigate the
stationarity of delay differential equations driven by L\'evy
processes of the form
\begin{equation}\label{eq.SDDEintro}
dX(t)=\left(\int_{[-\alpha,0]}X(t+s)\, \mu(ds)\right)\,dt +
F(X)(t-)\,dL(t).
\end{equation}
$L$ denotes a general L\'evy process, the drift term is obtained by integrating past
values with respect to a signed measure $\mu$ and  the nonlinear coefficient $F(X)$
depends on $(X(s): s\in [t-\alpha,t])$ at time $t$, see Section \ref{se.preliminaries}
for details. We do not consider a nonlinearity in the drift in order to concentrate on
the effects of the nonlinear noise term, which is facilitated by a variation of constants
formula. While the solution processes are not Markovian anymore, one can retrieve the
Markov property by regarding segments of the trajectories as processes in a function
space. The delayed noise term causes a fundamental degeneration of the segment process:
we show that the Markov semigroup is not Feller and not eventually strong Feller, but
eventually Feller. Consequently the uniqueness of an invariant measure can not be derived
by the strong Feller property.

Stationarity results for L\'evy-driven stochastic differential
equations, even in the non-delay case, are not so widespread.
Non-Gaussian stationary Ornstein-Uhlenbeck processes have been
attracting increasing attention recently due to their use in
financial modelling and the relationship with self-decomposable
distributions, cf. Barndorff-Nielsen and Shephard \cite{BNSh}.
Invariant measures for ordinary differential equations with a
nonlinear drift term and additive stable noise have been studied
analytically by Albeverio, R\"udiger and Wu \cite{ARW}, but our
general results, even when specified to the non-delayed case, seem
to be new.

The question of the existence of stationary solutions of
stochastic equations with delay goes back to the 60s in the work
of It\^o and Nisio \cite{ItNi64}. They have proved the existence,
but not the uniqueness of a stationary solution for Wiener-driven
delay differential equations under the condition that the drift is
obtained by a delayed perturbation of a stable instantaneous
feedback. For a more general non-linear drift functional and
additive white noise, Scheutzow \cite{Sch84} derived sufficient
conditions for the existence of an invariant probability measure
in terms of Lyapunov functionals. For a similar approach and
connections to stochastic partial differential equations see
Bakhtin and Mattingly \cite{BakMat03}. For the L\'evy-driven
equation \eqref{eq.SDDEintro} with constant $F$, that is additive
noise, Gushchin and K\"uchler \cite{GuKu98} have established the
existence and uniqueness of stationary solutions.

Our work is based on analyzing the segment process in a function
space and it is therefore closely related to results for
stochastic evolution equations in infinite-dimensional spaces. In
case of additive noise an extensive literature for the
stationarity of solutions of stochastic evolution equations
exists, see Da Prato and Zabczyk \cite{DaPZab96}. Much less is
known for non-additive noise, see for example Chow and Khasminskii
\cite{ChKhas97} for some general results. An infinite-dimensional
analogue of equation \eqref{eq.SDDEintro} driven by a Wiener
process is considered by Bonaccorsi and Tessitore \cite{BoTe01}.
They obtain a stationarity result for small Lipschitz constants by
a fixed point argument.

To prove the existence of a stationary solution of \eqref{eq.SDDEintro} under
rather
general conditions, we consider the segment process with values in the Skorokhod
space
$D([-\alpha,0])$. First, we establish the Feller property for the Markov
semigroup after
time $t=\alpha$. Under the main assumption of a stable drift, we establish the
tightness
of the solution segments using semimartingale characteristics and apply the
Krylov-Bogoliubov method to obtain an invariant measure on the Skorokhod space.
Due to
the absence of the strong Feller property Doob's method fails to prove
uniqueness of the
invariant measure. From an abstract point of view the loss of the strong Feller
property
is due to the degeneracy of the diffusion term when the equation is lifted to
the segment
space: the driving process is only one-dimensional, cf. Gatarek and Goldys
\cite{GataGold} for the abstract non-degenerate case. The question of uniqueness
of the
stationary solution turns out to be subtle and the degeneracy of the noise
process does
not permit a straight-forward analytical treatment. While for certain cases
uniqueness
will be shown to hold, a counterexample lets us suspect that uniqueness fails in
greater
generality. Nevertheless, the correlation structure of the solution process, if
it
exists, is uniquely determined and analytically tractable.

In the next section we briefly review some basic facts about stochastic delay
differential equations. Section 3 is devoted to the variation of constants
formula and
properties of the Markov semigroup. The existence and uniqueness of stationary
solutions
are discussed in Sections 4 and 5, respectively.

\section{Preliminaries} \label{se.preliminaries}

We follow standard notation, in particular we write $C[a,b]$ for the space of real-valued
continuous functions on $[a,b]$.
The Skorokhod space of all real-valued functions
on $[a,b]$ that are right-continuous and have left limits at every point ({\em
c\`adl\`ag} for short) is denoted by $D[a,b]$. It is endowed with the {\em Skorokhod
metric} $d_S$ given by
\begin{align*}
d_S(\phi,\psi):=\inf_{\lambda\in\Lambda[a,b]} \Big( \|\phi\circ \lambda -\psi
\|_\infty +
\|\lambda -\Id\|_\infty\Big),
\end{align*}
where
$
\Lambda[a,b]:= \{ \lambda:[a,b]\to [a,b]: \lambda\mbox{ is an increasing
homeomorphism}\}.
$
Note that $d_S(\phi_n,\phi)\to 0$ implies the convergences $\phi_n(a)\to
\phi(a)$,
$\phi_n(b)\to\phi(b)$, but not the pointwise convergence in the interior
$(a,b)$.

The space $(D[a,b],d_S)$ is a separable metric space. Moreover, there exists an
equivalent metric $d$ on $D[a,b]$ such that $(D[a,b],d)$ is a complete separable
metric
space, see for instance Jacod and Shiryaev \cite{JacShi03}. We endow $D[a,b]$
with the
corresponding Borel $\sigma$-algebra $\Borel(D[a,b])$. For $\phi\in D[a,b]$ we
denote by
$\phi(t-)$ its left-hand limit at $t$ and we define $\Delta\phi(t):=\phi(t)-
\phi(t-)$,
$t\in (a,b]$, and $\Delta\phi(a)=0$. For $\alpha>0$ and a function $\phi\in
D[-\alpha,\infty)$ we introduce the \textit{segment} of $\phi$ at time $t\ge 0$
as the
function
\begin{align*}
  \phi_{t}:[-\alpha,0]\to\R, \qquad \phi_{t}(u):=\phi(t+u).
\end{align*}

Let us first turn our attention to the deterministic delay equation underlying
the
stochastic equation \eqref{eq.SDDEintro}:
\begin{align}
  \begin{split}
     x(t)&= \phi(0)+\int_0^t \left(\int_{[-\alpha,0]} x(s+u) \, \mu(du)\right)ds
       \quad\text{for }t\ge 0,\\
     x(u)&= \phi(u)\quad \text{for }u\in [-\alpha,0],
  \end{split}\label{eq.det}
\end{align}
where $\mu$ is a signed finite Borel measure and the initial function $\phi$ is
in
$D[-\alpha,0]$. Note that the inner integral exists because $\phi$ and a
fortiori also
$x$ are measurable and locally bounded.

As the fundamental system in linear ordinary differential equations and the
Green
function in partial differential equations, the so-called \textit{fundamental
solution}
or \textit{resolvent} plays a major role in the analysis of \eqref{eq.det}. It
is the
function $r:\R\to \R$ which satisfies \eqref{eq.det} with the initial condition
$r(0)=1$
and $r(u)=0$ for $u\in [-\infty,0)$. The solution $x(\cdot,\phi)$ of
\eqref{eq.det} for
an arbitrary initial segment $\phi\in D[-\alpha,0]$ exists, is unique, and can
be
represented as
\begin{equation}\label{eq.xdet}
x(t,\phi)=\phi(0)r(t)+\int_{[-\alpha,0]} \int_{s}^0 r(t+s-u)\phi(u)\,du\,\mu(ds)
\quad\text{for }t\ge 0,
\end{equation}
cf. Chapter I in  Diekmann et al. \cite{DGVW95}. The fundamental solution
converges for
$t\to\infty$ to zero if and only if
\begin{align}\label{eq.v0}
 v_0(\mu)
 :=\sup\left\{\Re(\lambda):\, \lambda\in\C,\;\lambda - \int_{[-\alpha,0]}e^{\lambda s} \,
\mu(ds)=0\right\}
 <0,
\end{align}
where $\Re(z)$ denotes the real part of a complex number $z$. In this case the
decay is
exponentially fast and the zero solution of \eqref{eq.det} is uniformly
asymptotically
stable.

Let us fix a complete probability space $(\Omega,\F,P)$ with a
filtration $(\F_t)_{t\ge 0}$ satisfying the usual conditions. We
study the following stochastic differential equation with time
delay:
\begin{align}
\begin{split}
 dX(t)&= \left(\int_{[-\alpha,0]}X(t+s)\, \mu(ds)\right)\,dt
          + F(X)(t-)\,dL(t)\quad\text{for }t\ge 0,\\
  X(u)&= \Phi(u)\quad\text{for } u\in [-\alpha,0],
 \end{split}\label{eq.stoch}
\end{align}
where $\mu$ is a signed finite Borel measure and the initial process $(\Phi(u):\,u\in
[-\alpha,0])$ is assumed to have trajectories in $D[-\alpha,0]$ and to be
$\F_0$-measurable. The driving process $L=(L(t):\, t\ge 0)$ is a L\'evy process. We
denote its L\'evy-Khintchine characteristic by $(b,\sigma^2, \nu)$ with respect to the
truncation function $x\mapsto x\1_{[-1,1]}(x)$.

Turning to the specification of the nonlinear mapping $F$, we remark
that results for the existence and uniqueness of strong or weak
solutions of stochastic delay differential equations driven by
Brownian motion appear in different generalities: Mohammed
\cite{Mohammed84} provides a result under random functional
Lipschitz conditions, Mao \cite{Mao97} discusses in addition the
\textit{method of steps}, which provides a unique solution without a
regular dependence of the coefficients on values in the past,
Liptser and Shiryaev \cite{LipShi01} give general results for weak
solutions and It{\^o} and Nisio \cite{ItNi64} consider the existence
of weak solutions for equations with finite and infinite delay.
Since our equations are driven by L\'evy processes and the most
general conditions are not our concern here, we follow Protter
\cite{Pro90} and merely assume that the deterministic functional
$F:D[-\alpha,\infty)\to D[-\alpha,\infty)$ is \textit{functional
Lipschitz and autonomous}, i.e. it is continuous with respect to the
Skorohod topology and it satisfies for all $\phi_i\in
D[-\alpha,\infty)$,
$i=1,2$:\\[5pt]
(a) there exists a constant $K>0$, independent of $\phi_i$ and $t$,
such that
\begin{align}     \label{eq.lipcond}
 |F(\phi_1)(t)-F(\phi_2)(t)|\le K \sup_{t-\alpha \le s\le
t}|\phi_1(s)- \phi_2(s)|\quad \text{ for all }t\ge 0;
\end{align}
(b) $ F(\phi_1(s+\cdot))(t)=F(\phi_1)(t+s)\quad\text{for all } t,\,s\ge 0.$\\[2pt]

Equivalently, setting $\tilde{F}(\phi|_{[-\alpha,0]}):=F(\phi)(0)$ the two
conditions can
be stated as $F(\phi)(t)=\tilde{F}(\phi_t)$ with a functional $\tilde{F}$ which
is
Lipschitz continuous on $D[-\alpha,0]$ equipped with the supremum norm.

We can rewrite the differential equation \eqref{eq.stoch} as the
integral equation
\begin{equation}\label{eq.prot}
X(t)=\Phi(0)+\int_0^t G_\Phi(X)(s)\,ds+\int_0^t H_\Phi(X)(s-)\,dL(s)
 \quad\text{for }t\ge 0,
\end{equation}
when introducing $G_\phi,\,H_\phi:\, D[0,\infty)\to D[0,\infty)$
for $s\ge 0$ and $\psi\in D[0,\infty)$ by (abusing notation slightly)
\begin{eqnarray}
&& \quad G_\phi(\psi)(s)=\int_{[-\alpha,0]}
\left(\phi_s(u)\1_{[-\alpha,-s)}(u)+\psi_s(u)\1_{[-
s,0]}(u)\right)\,\mu(du),\label{eq.defG}\\
&& \quad H_\phi(\psi)(s)=F\left(\phi\1_{[-\alpha,0)} +\psi\1_{[0,\infty)}\right)(s).
   \label{eq.defH}
\end{eqnarray}
For ${\cal F}_0$-measurable initial segments $\Phi$ the mappings $G_\Phi$ and $H_\Phi$
are functional Lipschitz in the definition of Protter \cite{Pro90} and we can invoke
Theorem V.7 in Protter \cite{Pro90} which ensures a unique strong solution of
\eqref{eq.stoch}. Recall that a strong solution of \eqref{eq.stoch} is an adapted,
stochastic process $X$ with c{\`a}dl{\`a}g paths satisfying \eqref{eq.prot}. The solution
is called unique if all solutions are indistinguishable. We denote the solution by
$(X(t):\, t\ge -\alpha)$ or $(X(t,\Phi):\,t\ge -\alpha)$.

\begin{examples}\mbox{}
{\rm
\begin{enumerate}
\item[{\rm (a)}] The no-delay case: if $\mu=b\delta_0$, a point-mass at zero,
and $F(\phi)(t)=f(\phi(t))$, $t\ge -\alpha$,  then the
equation reads
\[ dX(t)=bX(t)\,dt+f(X(t))\,dL(t)\quad\text{ for }t\ge 0.
\]
If $f$ is Lipschitz continuous, then $F$ is easily seen to be
functional Lipschitz and autonomous.
\item[{\rm (b)}] The point-delay case: suppose $\mu=\sum_{i=1}^n
b_i\delta_{\alpha_i}$ and
$F(\phi)(t)=f(\phi(t-\alpha_1),\ldots,\phi(t-\alpha_n))$, $t\ge 0$, and
$F(\phi)(u)=F(\phi)(0)$, $u\in [-\alpha,0]$, with $\alpha_i\in [-\alpha,0]$.
Then the
equation reads
\[ dX(t)=\sum_{i=1}^n b_iX(t-\alpha_i)\,dt+f(X(t-\alpha_1),\ldots,X(t-
\alpha_n))\,dL(t)\text{ for }t\ge 0
\]
and $F$ is again autonomous and functional Lipschitz if $f$ is Lip\-schitz in all its
arguments.
\item[{\rm (c)}] The distributed-delay case: for $\mu(ds)=b(s)\,ds$ and
$F(\phi)(t)=f(\int_{[-\alpha,0]} \phi(t+s)c(s)\,ds)$, $t\ge 0$, and
$F(\phi)(u)=F(\phi)(0)$, $u\in [-\alpha,0]$, we obtain for $t\ge 0$
\[
dX(t)=\int_{[-\alpha,0]}X(t+s)b(s)\,ds\,dt+f\left(\int_{[-
\alpha,0]}X(t+s)c(s)\,ds\right)dL(t).
\]
Again, we need $f$ to be Lipschitz in order to have $F$ functional Lipschitz and
autonomous.
\item[{\rm (d)}] Further examples and counterexamples: other useful path-dependent mappings like
$F(\phi)(t)=\sup_{u\in [t-\alpha,t]}\phi(u)$ and their combinations with
Lipschitz
functions are functional Lipschitz and autonomous. Beware, however, that not all
Lipschitz continuous functionals $\tilde{F}$ on $D[-\alpha,0]$ give rise to a
functional
$F: D[-\alpha,\infty)\to D[-\alpha,\infty)$, for instance the jump size
functional
$F(\phi)(t)=\Delta \phi(t)$ is not c\`adl\`ag for c\`adl\`ag functions $\phi$
with jumps.
It is interesting to note that all admissible linear functionals are given by
$\tilde{F}(\phi)=\int_{[-\alpha,0]}\phi(u)\,\rho(du)$ with $\rho$ ranging
through the
space of finite Borel measures, which follows from the result by Pestman
\cite{Pest95}
when excluding the part based on jump sizes.

\end{enumerate}
}
\end{examples}


\section{Properties of the solution}

\subsection{The variation of constants formula}

Many of our considerations will be based on a stochastic
convolution equation, the {\em variation of constants formula}.
This formula is easily derived if the driving process has bounded
second moments, but no longer for processes where an It\^o
isometry or inequality fails. We provide a proof separately in
Rei{\ss} et al. \cite{Onno05}.
\begin{theorem}
Let $F$ be functional Lipschitz. Then for a stochastic process $X=(X(t):\,t\ge -
\alpha)$
and initial condition $\Phi$ the following are equivalent:
\begin{enumerate}
\item[{\rm 1)}] $X$ is the unique solution of \eqref{eq.stoch} with $X_0=\Phi$;
\item[{\rm 2)}] $X$ obeys the variation of constants formula:
\begin{align}\label{eq.varcons}
 X(t)= \begin{cases}\displaystyle
  x(t,\Phi) + \int_0^t r(t-s)F(X)(s-)\,dL(s), &t\ge 0,\\
   \Phi(t), &t\in [-\alpha,0],
   \end{cases}
\end{align}
where $r$ is the fundamental solution of equation \eqref{eq.det}.
\end{enumerate}
\end{theorem}

\subsection{Measurability of the segment process}

Our further work will be strongly based on considering the segment
process $(X_t:\, t\ge 0)$ in $D[-\alpha,0]$ instead of the
real-valued process $(X(t):\,t\ge - \alpha)$. This approach is
natural because the segment process is Markovian and turns out to
be eventually Feller. These properties will pave the way for our
further analysis.

In the case that $L$ is a Brownian motion the segment process is immediately a
Feller
process on the path space $C([-\alpha,0])$, see Theorem III.3.1 in Mohammed
\cite{Mohammed84}, which is not true in our setting because of the discontinuity
of the
shift semigroup on $D[-\alpha,0]$.

The following two Lemmas establish certain measurability and continuity properties of the
segment process. For a continuous path space similar  properties have been studied in
Chapter 3.7 of Da Prato and Zabczyk \cite{DaPZab92} and Lemma II.2.1 in Mohammed
\cite{Mohammed84} in Chapter 3.7.

\begin{lemma}\label{lem.msblsegment}
Let $Y=(Y(t):\, t\in [a,b])$ be a progressively measurable, real-valued
stochastic
process with c{\`a}dl{\`a}g  paths. Then $Y$ is a $D[a,b]$-valued random
variable.
\end{lemma}
\begin{proof}
According to \cite[Thm. 14.5]{Bil68} the Borel $\sigma$-algebra of
$D[a,b]$ coincides with the cylindrical $\sigma$-algebra generated
by all point evaluations $\phi\mapsto \phi(c)$ for $\phi\in D[a,b]$
and arbitrary $c\in [a,b]$. Hence $Y$ is measurable as $Y(c)$ is
Borel measurable for every $c\in [a,b]$.
\end{proof}

\begin{lemma}\label{lem.stochcon}
Let $(Y(t):\, t\ge 0)$ be a stochastically continuous process with
c{\`a}dl{\`a}g paths.
Then the segment process $(Y_t:\, t\ge \alpha)$ in $D[-\alpha,0]$ is
stochastically
continuous as well. Moreover, there exists a jointly measurable modification of
$(Y_t:\,
t\ge \alpha)$.
\end{lemma}

\begin{proof}
For $h>0$ we define the homeomorphism $\lambda_h:[-\alpha,0]\to
[-\alpha,0]$ by $\lambda_h(s):= s-h $ for $s\in [-\alpha+2h, -h]$
and affine respectively on $[-\alpha,-\alpha+2h]$ such that
$\lambda_h(-\alpha)=-\alpha$ and on $[-h,0]$ such that
$\lambda_h(0)=0$. Then $\|\lambda_{h}-\Id\|_\infty\le h$ and
\begin{align*}
Y_{t+h}(\lambda_h(s))
 = \begin{cases}
      Y(s+t), &  s\in [-\alpha+2h,-h] ,\\
       Y(t+h+\lambda_h(s)), & s\in [-\alpha,-\alpha+2h)\cup (-h,0].
    \end{cases}
\end{align*}
Therefore, we obtain
\begin{align*}
d_S(Y_{t+h},Y_t) \le  \norm{Y_{t+h}\circ \lambda_h -Y_t}_\infty
+\|\lambda_{h}-\Id\|_\infty
\to   \abs{\Delta Y(t)}\qquad\text{ as }\,h\downarrow 0.
\end{align*}
Hence, $t\mapsto Y_t$ is right-continuous at $t_0$ if $Y$ is
continuous at $t_0$. Similarly, one can establish
$\overline{\lim}_{h\downarrow 0}d_S(Y_{t-h},Y_t)\le \abs{\Delta
Y(t-\alpha)}$. We conclude that $t\mapsto Y_t$ is stochastically
continuous at $t_0$ if $P(\Delta Y(t_0-\alpha)\not= 0)=P(\Delta
Y(t_0)\not= 0)=0$, which follows from the stochastic continuity of
$Y$.

Any stochastically continuous process with values in a Polish
space has a jointly measurable modification, which is proved
following \cite[Prop. 3.2]{DaPZab92}, but measuring the distance
with the metric of this space. This gives the final assertion.
\end{proof}

\subsection{The Feller property}\label{sse.feller}

Basic tools for deriving the existence and uniqueness of invariant
measures are the Feller and strong Feller property of the Markov
semigroup defined by the segment process. We establish here the
Markov property, the Feller property after time $\alpha$ and give
examples that the immediate Feller and the eventually strong
Feller property fail in general. For our purposes the ordinary
Markov property of the segment process will be sufficient, but the
strong Markov property can also be derived following the lines of
Chapter 9.2 in Da Prato and Zabczyk \cite{DaPZab92}.

\begin{proposition}\label{pro.markov}
Let $X$ be the unique solution of \eqref{eq.stoch}. Then the segment process
$(X_t:\,t\ge
0)$ is a Markov process on $D[-\alpha,0]$:
\begin{align*}
 P(X_t\in B\,|\, \F_s)= P(X_t\in B\,|\, X_s)\quad\text{$P$-a.s. }
\end{align*}
for all $\;t\ge s\ge  0$ and Borel sets $B\in \Borel(D[-\alpha,0])$.
\end{proposition}

\begin{proof}
We fix $u\ge 0$ and consider for $t\ge u$ the equation
\begin{align*}
 X^{u,\phi}(t)&=\phi(0)+\int_u^t \int_{[-\alpha,0]}X^{u,\phi}(s+v)\, \mu(dv)\,ds
  + \int_u^t F(X^{u,\phi})(s-)\, dL(s), \\
 X^{u,\phi}(m)&= \phi(m-u)\quad\text{for }m\in [u-\alpha,u] \text{ and }
 \phi\in D[-\alpha,0].
\end{align*}
We denote the unique strong solution by $(X^{u,\phi}(t):t\ge u-\alpha)$ and the
segment
process by $(X^{u,\phi}_t:t\ge u)$.

We define ${\cal G}_u:=\sigma(L(s)-L(u):\, s\ge u)$ which is independent of the
$\sigma$-algebra $\F_u$ from the given filtration. The solution $X^{u,\phi}(t)$
is ${\cal
G}_u$-measurable for every $t\ge u$ and, as in the proof of Lemma
\ref{lem.msblsegment},
it follows that the segment $X^{u,\phi}_t$ is ${\cal G}_u$-measurable as well.
The
uniqueness of the solution implies $X(s)=X^{u,X_u}(s)$ for every $s\ge u-\alpha$
and thus
$X_t=X_t^{u,X_u}$ for every $t\ge u$ with probability one. By construction (cf.
\cite{Pro90}) the solution process depends in a measurable way on the initial
condition
so that the function
\begin{align*}
 A: D[-\alpha,0]\times \Omega\to\R,\qquad
  A(\phi,\omega):=\1_B( X_t^{u,\phi}(\omega))
\end{align*}
is  measurable for every $B\in \Borel(D[-\alpha,0])$  and
independent of $\F_u$ for fixed $\phi$. An application of the
factorisation lemma  \cite[Prop. 1.12]{DaPZab96} yields $P$-almost
surely
\begin{align*}
 P(X_t\in B\,| \F_u)
  = \E[\1_B(X_t^{u,X_u})\,| \F_u]
  = \E[A(X_u,\cdot)\, |\, \F_u]
   = \E[A(\phi,\cdot)]|_{\phi=X_u},
\end{align*}
which ends the proof because the right-hand side is
$\sigma(X_u)$-measurable.
\end{proof}

Let $B_b(D[-\alpha,0])$  denote the space of all real-valued bounded Borel
functions on
$D[-\alpha,0]$ endowed with the supremum norm $\norm{\cdot}_\infty$ and
$C_b(D[-\alpha,0])$ its subspace of continuous functions. Due to Proposition
\ref{pro.markov} the operators
\begin{align*}
 P_{s,t}: B_b(D[-\alpha,0])\to B_b(D[-\alpha,0]), \qquad
      P_{s,t}f(\phi):=\E[f(X_t^{s,\varphi})]
\end{align*}
have the property that $P_{u,s}P_{s,t}=P_{u,t}$ for $0\le u\le
s\le t$. By homogeneity we have $P_{s,t}=P_{0,t-s}$ for $0\le s\le
t$, cf. Thm. V.32 in Protter \cite{Pro90}, and the operators
$P_t:=P_{0,t}$, $t\ge 0$, form a {\em Markovian semigroup}. The
Markovian semigroup will be called {\em eventually Feller} if
there exists a $t_0\ge 0$ such that for any $f\in
C_b(D[-\alpha,0])$ the following two conditions are satisfied:
\begin{align}
 &P_t f \in C_b(D[-\alpha,0]) \qquad\text{for every }t\ge t_0,
 \label{eq.feller}\\[2pt]
 &\lim_{s\downarrow t} P_sf(\phi)=P_tf(\phi)\qquad\text{for every }
  \phi \in D[-\alpha,0],\,t\ge t_0\label{eq.stochcont}.
\end{align}
By definition the solution $(X(t):\,t\ge 0)$ has c\`adl\`ag paths
and it is easily observed from (\ref{eq.prot}) that it has no
fixed times of discontinuity: $P(\Delta X(t)\not=0)\le P(\Delta
L(t)\not=0)=0$. Hence, the process $(X(t):\,t\ge 0)$ is
stochastically continuous. By Lemma \ref{lem.stochcon} so is the
segment process $(X_t:\,t\ge \alpha)$ and thus condition
\eqref{eq.stochcont} is fulfilled for $t_0=\alpha$. The semigroup
is not stochastically continuous for $t_0<\alpha$ and condition
\eqref{eq.feller} fails for $t_0<\alpha$ due to the discontinuity
of the shift semigroup, as the following example demonstrates.

Consider $0<\beta<\alpha$ and the initial functions $\phi^n:=\1_{[-\beta(1-n^{-
1}),0]}$
which for $n\to\infty$ converge in $D[-\alpha,0]$ to $\phi^\infty:=\1_{[-
\beta,0]}$. The
corresponding solution segments $X^n_t$, for an arbitrary specification of $F$
and $L$ in the
differential equation, satisfy $X^n_{\alpha-\beta}(-\alpha)=\phi^n(-\beta)=0$,
while
$X^\infty_{\alpha-\beta}(-\alpha)=1$ holds. Hence,
$d_S(X^n_{\alpha-\beta},X^\infty_{\alpha-\beta})\ge 1$, which implies that
$\phi\mapsto
P_tf(\phi)$ is not continuous for $f(\psi):=\abs{\psi(-\alpha)}\wedge 1\in
C_b(D[-\alpha,0])$ and any time $t\in (0,\alpha)$. Similarly, $t\mapsto
P_tf(\phi^\infty)$ is seen to be discontinuous at $t=\alpha-\beta$.

We now establish condition \eqref{eq.feller} for $t_0=\alpha$ by showing even
more,
namely that $\phi\mapsto (X(t,\phi):\,t\in [0,T])$ is continuous from $D[-
\alpha,0]$ to
the space of c\`adl\`ag processes with the uniform convergence on $[0,T]$ in
probability,
which is stronger than convergence in the Skorokhod topology in law. We start
with a norm
estimate in spirit of \'Emery's inequality  before proving the main result. In
accordance
with Section V.2 in Protter \cite{Pro90} we employ the following norms for
semimartingales $(Z(t):t\ge 0)$ and adapted c\`adl\`ag processes $(Y(t):t\ge
0)$:
\begin{align*}
&\norm{Y}_{S^2[0,T]}^2:=\E\big[\sup_{0\le t\le T}Y(t)^2\big],\\
&\norm{Z}_{H^2[0,T]}^2:=\inf\left\{\E[M,M]_T+\E[\TV(A)(T)^2]\right\},
\end{align*}
where the infimum is taken over all possible decompositions $Z=M+A$ where $M$ is a local
martingale and $A$ a bounded variation process with $M(0)=A(0)=0$. The total variation of
$A$ on $[0,T]$ is denoted by $\TV(A)(T)$.

The quadratic variation process is defined by $[Z,Z]:=Z^2-\int Z(s-)\,dZ(s)$. Based on
these norms the spaces $H^2[0,T]$ and $S^2[0,T]$ are constructed canonically. Moreover,
they are Banach spaces, and $H^2[0,T]$ is continuously embedded in $S^2[0,T]$.
\begin{lemma}\label{lem.emery}
Suppose the L\'evy process $L$ has a finite second moment and
$(H(t):\,0\le t\le T)$ is an adapted c\`adl\`ag process with
$\int_0^T \E[H(t-)^2]\,dt<\infty$. Then
\begin{align*}
&\norm{\int_0^\cdot H(s-)\,dL(s)}_{H^2[0,T]}^2\\
&\qquad\qquad\le \left(\sigma^2+\int x^2\,\nu(dx)+(\E L(1))^2T\right)\int_0^T
\E[H(t-)^2]\,dt.
\end{align*}
\end{lemma}

\begin{proof}
This follows from the decomposition $L(t)=M(t)+t \E L(1)$ with $M$ a square integrable
martingale.
\end{proof}

The surprising result of the next proposition, which says that convergence of
the initial conditions
in Skorokhod metric implies uniform convergence of the solution processes, is
essentially
due to the fact that the driving L\'evy process is a semimartingale without
fixed time of
discontinuity.

\begin{proposition}\label{pro.feller}
Assume $F:D[-\alpha,\infty)\to D[-\alpha,\infty)$ is continuous with respect to
the
Skorokhod metric. Let $X^n$ be the solution of equation \eqref{eq.stoch} with
deterministic initial segment $\phi^n$, let $\phi^n\to \phi$ in $D[-\alpha,0]$
and let
$X$ be the solution with initial segment $\phi$. Then $(X^n(t):t\ge 0)$
converges to
$(X(t): t\ge 0)$ uniformly on compact sets in probability.
\end{proposition}

\begin{proof}
We consider first a stopping time $R$ such that the process $L^{R-}$ is
$\alpha$-sliceable for
some suitably small $\alpha>0$ in the sense of \cite{Pro90}.

In analogy to \cite[Thm. V.10]{Pro90} we use the representation
\eqref{eq.prot} and put
\begin{align*}
Y^n(t)&:=\int_0^t (G_\phi(X)-G_{\phi^n}(X))(s)\,ds\\
  &\qquad +\int_0^t (H_\phi(X)-H_{\phi^n}(X))(s-)\,dL^{R-}(s), \\
{\cal G}^n(U)(t)&:=G_{\phi^n}(X)(t)-G_{\phi^n}(X-U)(t),\\
{\cal H}^n(U)(t)&:=H_{\phi^n}(X)(t)-H_{\phi^n}(X-U)(t),
\end{align*}
$t \ge 0$, to obtain for $U^n:=X-X^n$ the equation
\[ U^n(t)=\phi(0)-\phi^n(0)+Y^n(t)+\int_0^t {\cal
G}^n(U^n)(s)\,ds+\int_0^t {\cal H}^n(U^n)(s-)\,dL^{R-}(s).
\]
By \cite[Lemma V.3.2]{Pro90}, extended to two driving semimartingales, the
solution $U^n$
of this equation satisfies $\norm{U^n}_{S^2[0,T]}\le
C\norm{\phi(0)-\phi^n(0)+Y^n}_{S^2[0,T]}$ for any $T>0$ with a constant $C>0$
depending
on the process $L^{R-}$ and a uniform bound for the Lipschitz constants of
$G_{\phi^n}$
and $H_{\phi^n}$. The Skorokhod metric ensures $\phi^n(0)\to \phi(0)$, so that
$\norm{U^n}_{S^2[0,T]}\to 0$ follows if $Y^n$ tends to zero in $S^2[0,T]$. The
latter is fulfilled if
\begin{equation}\label{eq.help1}
\E\left[\int_0^T (G_\phi(X)(t)-G_{\phi^n}(X)(t))^2+ (H_\phi(X)(t)-
H_{\phi^n}(X)(t))^2
dt\right]
\end{equation}
tends to 0 as $n\to\infty$, due to the continuous embedding $H^2\hookrightarrow S^2$ from
\cite[Thm.V.2]{Pro90} and Lemma \ref{lem.emery} with the additional observation that
$\|\int_0^\cdot J(s-)dL^{R- }(s)\|_{H^2[0,T]}\le \|\int_0^\cdot J(s-)dL(s)\|_{H^2[0,T]}$
for any process $J\in S^2[0,T]$. Let $\omega$ be fixed for the moment. The functions
$H_{\phi^n}(X(\omega))$ converge in the Skorokhod topology to $H_{\phi}(X(\omega))$,
which implies convergence in $L^2[0,T]$. Concerning the sequence $G_{\phi_n}$ we have
\begin{align*}
&\int_0^T (G_{\phi_n}(X(\omega))(t)-G_{\phi}(X(\omega))(t))^2\,dt
 =\int_{[-\alpha,0)}\int_{[-\alpha,0)} \int_0^T (\phi_n(t+u)\\
 &\quad-\phi(t+u)) (\phi_n(t+v)-\phi(t+v))\1_{[-\alpha,-t)}(u)
 \1_{[-\alpha,-t)}(v)\,dt\,\mu(du)\,\mu(dv).
\end{align*}
This expression converges to zero as $n\to\infty$, since the
Skorohod convergence of $\phi_n$ to $\phi$ implies $\phi_n\to\phi$
Lebesgue a.e.\ and the sequence $(\phi_n)_n$ is uniformly bounded.
Again by \cite[Lemma V.3.2]{Pro90} the solution process $X$ is an
element of $S^2[0,T]$, whence by the uniform linear growth of
$(G_{\phi^n})$ and $(H_{\phi^n})$ the argument inside the
expectation in \eqref{eq.help1} is dominated by a $P$-integrable
function. The Dominated Convergence Theorem thus gives the
convergence in \eqref{eq.help1} such that
$\norm{X-X^n}_{S^2[0,T]}=\norm{U^n}_{S^2[0,T]}\to 0$ for any $T>0$.

Next, let $L$ be an arbitrary L\'evy process. According to
\cite[Theorem V.5, p.192]{Pro90} there exist stopping times $0=
T_0\le T_1\le T_2\le\cdots$ such that $\sup_\ell T_\ell=\infty$
a.s.\ and $L^{T_\ell-}$ is $\alpha$-sliceable for each $\ell$.
Consider equation (\ref{eq.stoch}) with $L$ replaced by
$L^{T_\ell-}$ and let $X^{n,\ell}$ denote the solution with
initial segment $\varphi_n$ and let $X^{\infty,\ell}$ denote the
solution with intial segment $\varphi$, for $n,\ell\in\N$. We have
shown above that $X^{n,\ell}\to X^{\infty,\ell}$ uniformly on
compact sets in probability for every $\ell$. Further, it is clear
from the equation that $X^{n,\ell}=(X^n)^{T_\ell-}$. Let now
$t>0$, $r>0$, and $\varepsilon>0$ be arbitrary. Choose an $\ell$
such that $P(T_\ell<t)<\varepsilon/2$. Then
\begin{align*}
&P(\sup_{0\le s\le t} | X^n(s)-X(s)|\ge r)\\
&\qquad\le P(\sup_{0\le s\le t} |X^n(s)-X(s)|\ge r\mbox{ and }T_\ell\ge
t)+P(T_\ell<t)\\
&\qquad\le P(\sup_{0\le s\le t} | X^{n,\ell}(s)-X^{\infty,\ell}(s)|\ge
r)+\varepsilon/2<\varepsilon
\end{align*}
for $n$ large. Hence $X^n\to X$ uniformly on compact sets in probability.
\end{proof}

Let us finally show that in general we cannot expect that the solution is
\textit{eventually strongly Feller}, which is characterised by the existence of
a $t_0>0$
such that for all $f\in B_b(D[-\alpha,0])$
\begin{align*}
 P_t f\in C_b(D[-\alpha,0])\qquad\text{for every }t\ge t_0.
\end{align*}
Using indicator functions for $f$, this implies
\begin{align*}
 \phi\mapsto P(X_t(\phi)\in B)\in C_b(D[-\alpha,0])
 \qquad\text{for every }B\in \Borel(D[-\alpha,0]),\,t\ge t_0.
\end{align*}

Suppose the functional $F$ in the equation \eqref{eq.stoch} is of the form
$F(\psi)(t)=f(\psi(t-\alpha))$ for $t\ge 0$ and $F(\psi)(t)=0$ for $t<0$, with a
Lipschitz-continuous homeomorphism $f:\R\to (a,b)$,
$b>a>0$, and consider the case that $L$ is standard Brownian motion. Then the
quadratic
variation $\langle X_t\rangle $ of the solution segment $X_t$, $t\ge \alpha$,
satisfies
\[ \langle X_t\rangle_u=\int_0^{t+u}
f^2(X(s-\alpha))\,ds \quad \text{$P$-a.s. for }u\in [-\alpha,0].
\]
Since both sides of the equation are continuous in $u$ for continuous $X$, there
is one
$P$-null exception set for all $u\in [-\alpha,0]$.
Consider the map
\begin{align*}
V(\phi)(u):=f^{-1}\left(\left(\frac{d\langle
\phi\rangle_u}{du}\right)^{1/2}\right),
\quad u\in [-\alpha,0],
\end{align*}
defined on the functions $\varphi$ with finite quadratic variation such that
$\frac{d\langle \phi\rangle_u}{du}\in (a^2,b^2)$ for Lebesgue-almost every $u\in
[-\alpha,0]$.
We have $P(V(X_t)(u)=X(t+u-\alpha),\,u\in [-\alpha,0])=1$ for all $t\ge \alpha$.
Iterating this map, we can recover with probability one the initial segment
$X_0$ from
observing $X_{m\alpha}$ since $V^m(X_{m\alpha})=X_0$ for every integer $m$. This
identifiability property shows that the laws of the segments
$X_{m\alpha}(\phi_1)$ and
$X_{m\alpha}(\phi_2)$ for different initial segments $\phi_1$ and $\phi_2$ must
be
singular. Hence, there is a contradiction to the strong Feller property at
$t_0=m\alpha$,
which asserts the continuous dependence of the laws on the initial condition.

In fact, this example even shows that the Markov semigroup is not
eventually regular in the sense of Da Prato and Zabczyk
\cite{DaPZab96}. We shall see in Section \ref{sec.nondelaydiff}
that this counterexample is due to the delay in the diffusion
coefficient.

\section{Existence of a stationary solution}\label{sec.statexist}

\subsection{Tightness}

We establish the tightness of the laws $\{\law(X_t)\}_{t\ge 0}$ in
${D}[0,\alpha]$ by considering the semimartingale characteristics.
Recall that $(b,\sigma^2, \nu)$ denotes the L\'evy-Khintchine
characteristic of  the L\'evy process $L$.

\begin{assumption}\label{AssStat}
\mbox{}
\begin{enumerate}
\item[{\rm (a)}] The delay measure $\mu$ in the drift satisfies $v_0(\mu)<0$
with
$v_0$ from equation \eqref{eq.v0}.
\item[{\rm (b)}] The jump measure $\nu$ satisfies
$\int_{\abs{x}>1} \log\abs{x}\,\nu(x)<\infty$.
\item[{\rm (c)}] The coefficient $F$ in equation \eqref{eq.stoch} is
functional Lipschitz, uniformly bounded and autonomous.
\end{enumerate}
\end{assumption}

Condition (a) yields the exponential decay of the fundamental solution, while
condition (b) ensures that $\int_0^t f(s)\,dL(s)$, for exponentially decaying
functions $f$ of locally bounded variation, converges in law and is already for
constant $F$ necessary for the existence of a stationary solution, as was shown
by Gushchin and K\"uchler \cite{GuKu98}, cf. also
Thm. 4.3.17 in Applebaum \cite{App04}. In condition (c)  restrictions on $F$
are imposed such that the differential equation is autonomous, has a unique
solution and the impact of the driving process cannot become too large. For the
latter the imposed boundedness of $F$ can certainly be relaxed considerably,
but will then depend on the large jumps of $L$, that is, on fine properties of
$\nu$.

\begin{proposition}\label{PropLevyMarg}
Grant Assumption \ref{AssStat}. Then the solution process
$(X(t):\,t\ge -\alpha)$ of \eqref{eq.stoch} with initial condition
$X_0=0$ has one-dimensional marginal laws $\{\law(X(t))\}_{t\ge
0}$ that are tight.
\end{proposition}

\begin{proof}
Let us split the L\'evy-process $L$ into two
parts, one of them consisting of jumps of size larger than one:
\begin{align*}
L(t)=N(t)+R(t)\qquad\text{with } N(t)=\sum_{s\le t} \Delta L(s)\1_{\{\abs{\Delta
L(s)}>1\}}.
\end{align*}
Then the variation of constants  formula \eqref{eq.varcons} yields $X=Y+Z$ with
\begin{align*}
Y(t)&:=\int_0^t r(t-s)F(X)(s-)\,dN(s) ,\quad t\ge 0, \\
\text{and} \qquad Z(t)&:=\int_0^t r(t-s)F(X)(s-)\,dR(s), \quad t\ge 0.
\end{align*}
Tightness of $(X(t):\, t\ge 0)$ will follow from tightness
of $Y$ and $Z$.

The fundamental solution $r$ decays exponentially with $\abs{r(t)}\le c e^{-
\beta t}$ for
some constants $c$, $\beta>0$ due to Assumption \ref{AssStat}(a). Considering
$Y$ first,
we obtain for any $K>0$ with $m:=\sup_{\psi}\abs{F(\psi)(0)}$ by time reversal
for the
compound Poisson process $N$ the estimate
\begin{align*}
 P(\abs{Y(t)}>K)
  & \le P \left(\sum_{s\le t} \abs{r(t-s)F(X)(s-)\Delta N(s)}>K \right)\\
  & \le P \left(\sum_{s\le t} c e^{-\beta(t-s)} m\abs{\Delta N(s)}>K\right)\\
  &= P \left (\sum_{s\le t} e^{-\beta s}\abs{\Delta N(s)}> \frac{K}{cm}\right).
\end{align*}
The tightness for $Y$ follows from the tightness of $\sum_{s\le t} e^{-\beta
s}\abs{\Delta N(s)}$ which has been established in \cite[Lemma 4.3]{GuKu98}
under
Assumption \ref{AssStat}(b).

Since $R$ is a L{\'evy} process with bounded jumps its canonical
decomposition is, by means of \cite[p.103]{App04}, of the simple
form $R(t)= R_0(t) + t\E{R(1)}$ where $(R_0(t):\, t\ge 0)$ is a
square-integrable martingale. We split $Z$ into the sum
$Z=Z_0+Z_1$ with
\begin{align*}
  Z_0(t)&:=\int_0^t r(t-s)F(X)(s-)\,dR_0(s),&& t\ge 0, \\
\text{and} \qquad
  Z_1(t)&:=\E[R(1)]\int_0^t r(t-s)F(X)(s-)\,ds, && t\ge 0.
\end{align*}
For $Z_1$ we easily obtain $P(\abs{Z_1(t)}>K) \le P(cm \beta^{-1}  |\E[R(1)]|>K)
=0$ for
$K$ sufficiently large, implying tightness of $Z_1$. As in Lemma \ref{lem.emery}
we
obtain
\begin{align*}
\E[Z_0(t)]^2\le
 \left(\sigma^2+\int_{\abs{x}\le 1} x^2\,\nu(dx)\right)m^2
      \int_0^t r^2(t-s)\,ds,\quad t\ge 0.
\end{align*}
Hence by the exponential decay of $r$, the sequence $(Z_0(t))_{t\ge 0}$ is
bounded in
$L^2_P(\Omega)$ and thus tight.
\end{proof}

\begin{proposition}\label{PropLevySeg}
In the setting of Proposition \ref{PropLevyMarg} we have that the laws
$\{\law(X(t+s)-X(t),\,s\in [0,\alpha])\}_{t\ge 0}$ are tight in $D[0,\alpha]$.
\end{proposition}

\begin{proof}
We are led to consider for $t\ge 0$ and $s\in [0,\alpha]$
\begin{align*}
 Y_t(s)&:=X(t+s)-X(t)\\
       &=\int_t^{t+s} \left(\int_{[-\alpha,0]} X(u+v)\,\mu(dv)\right)\,du
        +\int_t^{t+s} F(X)(u-)\,dL(u).
\end{align*}
Let us introduce for $t\ge 0$ the semimartingale $(I_t(s):\, s\in [0,\alpha])$
by letting
\begin{align*}
 I_t(s):= \int_t^{t+s} F(X)(u-)\,dL(u)\qquad\text{for }s\in [0,\alpha].
\end{align*}
Now, either by following the lines in \cite[III.2.c]{JacShi03} and using the
L\'evy-It\^o decomposition or by applying \cite[Prop. 7.6]{Rai00} the
semimartingale characteristic $(B_{I_t},C_{I_t},\nu_{I_t})$ of $(I_t(s):s\in
[0,\alpha])$ is found to be
\begin{align*}
B_{I_t}(s)&=\int_t^{t+s} \big(bF(X)(u-)\\
 &\quad  + \int xF(X)(u-)\left(\1_{(-1,1)}(xF(X)(u-))- \1_{(-
1,1)}(x)\right)\nu(dx)\big)du,\\
C_{I_t}(s)&=\sigma^2 \int_t^{t+s} F^2(X)(u-)\,du,\\
\nu_{I_t}(ds,dx)&=ds\times K_{I_t}(X,t+s, dx),
\end{align*}
where for $y\in D[-\alpha,\infty),\,u\ge t$ and a Borel set $A\in \Borel(\R)$
\[K_{I_t}(y,u,A):=\int \1_{A\setminus
\{0\}}\left(F(y)(u-)x\right) \,\nu(dx).
\]
Hence, the semimartingale $(Y_t(s):\, s\in [0,\alpha])$ has the
characteristic $(B_{Y_t},C_{Y_t},\nu_{Y_t})$ with
$C_{Y_t}=C_{I_t}$, $\nu_{Y_t}=\nu_{I_t}$ and
\begin{align*}
 B_{Y_t}(s)=B_{I_t}(s)+ \int_t^{t+s} \left(\int_{[-\alpha,0]}
X(u+v)\,\mu(dv)\right)du .
\end{align*}
We prove the tightness of $(Y_t)_{t\ge 0}$ by means of \cite[Thm.
VI.4.18]{JacShi03} and \cite[VI.4.20]{JacShi03}. For that we have to verify
that
\begin{align*}
 a_{Y_t}(s):= \TV(B_{Y_t})(s) + C_{Y_t}(s)
     + \int_{[0,s]\times \R} (\abs{x}^2\wedge 1)\,\nu_{Y_t}(du,dx),\quad s\in
[0,\alpha],
\end{align*}
forms a tight sequence $(a_{Y_t})_{t\ge 0}$ of processes  and all limit points
of the
sequence $\{\law(a_{Y_t})\}_{t\ge 0}$ as $t\to\infty$ are laws of continuous
processes.
According to \cite[Prop.\ VI.3.33 and VI.3.35]{JacShi03} this will follow if
there exist
some increasing processes $(A_{Y_t})_{t\ge 0}$ satisfying these conditions and
in
addition $A_{Y_t} -a_{Y_t}$ defines for every $t\ge 0$ an increasing process,
since
$a_{Y_t}(s)\ge 0$ a.s.\ for all $s\in [0,\alpha]$, $t\ge 0$. To obtain such
processes
$(A_{Y_t}(s):\, s\in [0,\alpha])$, we estimate
\begin{align*}
\TV(B_{Y_t})(s) &= \int_t^{t+s}\bigg|\int_{[-\alpha,0]} X(u+v)\,\mu(dv) + b
F(X)(u-)
 \\
&  + \int_{\R} xF(X)(u-)  \left(\1_{(-1,1)}\left(xF(X)(u-)\right)-
      \1_{(-1,1)}(x)\right)\nu(dx)\bigg| du\\
&\le \int_t^{t+s}\left( \abs{\int_{[-\alpha,0]} X(u+v)\,\mu(dv)} + \abs{b}m +
c\right)\,du ,
\intertext{where $m:=\sup_{\psi}\abs{F(\psi)(0)}$ and the finite
constant $c$ is defined by} &c:= \int_{\frac{1}{m}\le \abs{x}<1}
  m\abs{x}\,\nu(dx)+\nu(\R\setminus (-1,1)).
\end{align*}
Therefore, the process $a_{Y_t}$ is majorized by
\begin{align*}
A_{Y_t}(s):&= \int_t^{t+s} \abs{\int_{[-\alpha,0]} X(u+v)\,\mu(dv)}\,du  \\
 &\qquad + s\left( \abs{b}m + c+ \sigma^2 m^2 + \int_{\R} ((m^2 x^2)\wedge 1)\,
\nu(dx)\right)\
 \;\text{for }s\in [0,\alpha]
\end{align*}
and $A_{Y_t}-a_{Y_t}$ is increasing. Since $A_{Y_t}$ is continuous and only
depends on
$t$ in the first term, it suffices to prove tightness in $C[0,\alpha]$ of the
first term:
\begin{align*}
 J_{Y_t}(s):=\int_t^{t+s}\babs{\int_{[-\alpha,0]}
X(u+v)\,\mu(dv)}\,du\quad\text{for }
  s\in [0,\alpha].
\end{align*}
Recalling $r(u)=0$ for $u<0$, we obtain by the variation of
constants formula
\begin{align*}
I(u)&:=\int_{[-\alpha,0]} X(u+v)\,\mu(dv)\\
&=\int_{[-\alpha,0]}\left( \int_0^u r(u+v-s)F(X)(s-)\,L(ds)\right)\,\mu(dv)\\
&=\int_0^u \dot{r}(u-s)F(X)(s-)\,dL(s).
\end{align*}
To prove tightness of the absolutely continuous processes
$(J_{Y_t})_{t\ge 0}$, it suffices to show that the process $(I(t+s): s\in
[0,\alpha])_{t\ge 0}$
is bounded in probability in $C[0,\alpha]$, that is,
\begin{align}\label{eq.conditItight}
\lim_{K\to\infty} \sup_{t\ge 0}
    P\left(\sup_{t\le u\le t+\alpha}\abs{I(u)}>K\right)= 0 .
\end{align}
Note that we have the exponential decay estimate
$\abs{\dot{r}(t)}\le c^\prime e^{-\beta t}$. Decomposing $L$ into
its drift, diffusion, and large and small jump parts, it is clear that
only integration with respect to the large jump part $N$ may pose
problems. As in the proof of Proposition \ref{PropLevyMarg},
however, the restriction on the large jumps in $L$ and the finite
intensity of $N$ yield in a similar manner
\begin{align*}
\sup_{t\le u\le t+\alpha}\babs{\int_0^u \dot{r}(u-s)F(X)(s-)\,dN(s)}
 &\le \sum_{s\le t+\alpha} c^\prime m e^{-\beta (t-s)}\abs{\Delta N(s)}
\end{align*}
and the tightness of the right-hand side by \cite[Lemma 4.3]{GuKu98}.
Thus we infer (\ref{eq.conditItight}).
\end{proof}

\begin{theorem}\label{th.tightness}
Grant Assumption \ref{AssStat}. Then for the solution process $(X(t):\, t\ge
-\alpha)$ of \eqref{eq.stoch} with initial condition $X_0=0$ the laws of the
segments $\{\law(X_t)\}_{t\ge\alpha}$ are tight in $D[-\alpha,0]$.
\end{theorem}

\begin{proof}
If we let $Z_t(s):=X(t-\alpha)$ for $s\in [-\alpha,0]$, then the processes
$(Z_t)_{t\ge
0}$ of constant functions are tight in $C[-\alpha,0]$ by Proposition
\ref{PropLevyMarg}.
On the other hand, $\{\law(X_t-Z_t)\}_{t\ge\alpha}$ are tight in $D[-\alpha,0]$
by
Proposition \ref{PropLevySeg} applying the time shift $t\mapsto t-\alpha$.
Therefore the
sum $(X_t-Z_t)+Z_t$ is tight in $D[-\alpha,0]$ using the result in
\cite[VI.3.33(a)]{JacShi03}.
\end{proof}

\subsection{From tight solutions to stationary solutions}

We use the construction due to Krylov and Bogoliubov, see for example Da Prato and
Zabczyk \cite{DaPZab96}. For the reader's convenience we include a complete proof, which
is tailored for our purposes. Consider equation (\ref{eq.stoch}) and its Markovian
semigroup $(P_t)_{t\ge 0}$ as defined below Proposition~\ref{pro.markov}. Denote by
$\calP=\calP(D[- \alpha,0])$ the set of Borel probability measures on $D[-\alpha,0]$,
endowed with the topology of weak convergence of measures. Let
$\langle\cdot,\cdot\rangle$ denote the duality pairing of $\calP$ and
$B_b:=B_b(D[-\alpha,0])$ given by $\langle \zeta,f\rangle = \int fd\zeta$,
$\zeta\in\calP$, $f\in B_b$. Define for $t\ge 0$ and $\zeta\in\calP$ the functional
$P_t^\ast\zeta$ by
$$
(P_t^\ast\zeta)f:= \langle  \zeta,P_t f\rangle,\quad f\in B_b.
$$
If $\zeta$ is the distribution of an initial segment $\Phi$, then $P_t^\ast\zeta$ is the
distribution of $X_t(\Phi)$, since
$$
\langle P_t^\ast\zeta,f\rangle =\int \E[f(X_t(\varphi))]\,\zeta(d\varphi)
 =\E\left[\E[f(X_t(\Phi))|\F_0]\right]=\E[f(X_t(\Phi))],
$$
for $f\in B_b$. A measure $\zeta\in\calP$ is called an \emph{invariant measure} or
\emph{stationary distribution} of (\ref{eq.stoch}) if $P_t^\ast\zeta=\zeta$ for all $t\ge
0$, that is, $\langle \zeta,P_t f\rangle=\langle \zeta,f\rangle$ for all $f\in B_b$ and
all $t\ge 0$.

It follows from Lemma~\ref{lem.stochcon} that $t\mapsto P_t^\ast \zeta$ is a continuous
map from $[\alpha,\infty)$ to $\calP$ and moreover $P_{s+t}^\ast\zeta=P_s^\ast
P_t^\ast\zeta$ for $s,t\ge 0$. Further, $P_t$ maps $C_b:=C_b(D[-\alpha,0])$ into $C_b$
for all $t\ge \alpha$, by Proposition\ \ref{pro.feller}.

Because of Theorem~\ref{th.tightness}, the next theorem follows from
Theorem~\ref{th.tighttoexist} below.

\begin{theorem}\label{th.statexists}
Grant Assumption~\ref{AssStat}. Then there exists a stationary distribution for
\emph{(\ref{eq.stoch})}.
\end{theorem}

\begin{theorem}\label{th.tighttoexist}
If for some $\zeta\in\calP$ the set $\{P_t^\ast\zeta:t\ge \alpha\}$ is tight, then there
exists an $\eta\in\calP$ such that $P_t^\ast \eta=\eta$ for all $t\ge 0$. Moreover,
$\eta$ is an element of the closed convex hull of $\{P_t^\ast\zeta:t\ge \alpha\}$ in
$\calP$.
\end{theorem}

\begin{proof}
Denote for convenience $T_t:=P_{t+\alpha}$ and $\zeta(t):=T_t^\ast\eta$, $t\ge 0$. First
we show that for each $t\ge 0$ there exists a unique $\theta_t\in\calP$ such that
$$
\langle\theta_t,f\rangle= \frac{1}{t} \int_0^t \langle \zeta(s),f\rangle ds
\quad\mbox{for all }f\in C_b.
$$
It is routine to show the uniqueness. In order to show existence, define
\begin{align*}
\phi(f):= \frac{1}{t}\int_0^t \scapro{\zeta(s)}{f}\, ds
  \qquad\text{for every } f\in C_b.
\end{align*}
Then $\phi\in C_b^\ast$, where $C_b^\ast$ denotes the dual Banach space of $C_b$. Let
$\epsilon>0$ and take $H\subseteq D[-\alpha,0]$ compact such that $\zeta(s)(H)\ge
1-\epsilon$ for all $s\in [0,t]$. For $f\in C_b$ with $\norm{f}_{\infty} \le 1$ and
$f(x)=0$ for all $x\in H$ we then have
\begin{align*}
\abs{\scapro{\zeta(s)}{f}} = \abs{ \int f \, d\zeta(s)} \le
  \norm{f}_\infty\, \zeta(s)(D[-\alpha,0]\setminus H)\le \epsilon
\end{align*}
for all $s\in [0,t]$, so $\abs{\phi(f)}\le \epsilon$. Now by the Riesz-Bourbaki
representation theorem (see for instance \cite[Prop. 5.2.5]{Bourbaki} and \cite[Prop.
5.6.12]{Bourbaki}), there exists a tight finite positive Borel measure $\theta$ on
$D[-\alpha,0]$ such that $\phi(f)= \scapro{\theta}{f}$ for every $f\in C_b$. Notice that
$\scapro{\theta}{\1_{D[-\alpha,0]}} =1$, so $\theta\in \calP$.

Next we show that $\theta$ is an element of the closure of the convex hull of
$\{\zeta(s):0\le s\le t\}$ in $\calP$. Let $M$ denote the weak* closure of the convex
hull of $\{\zeta(s): 0\le s\le t\}$ in $C_b^\ast$. Then $M$ is a weak* closed convex set
in $C_b^\ast$. The Hahn-Banach Theorem
implies that for any $\psi\in C_b^\ast\setminus M$ there exist $f\in C_b$ and $\beta\in
\R$ such that $\psi(f)<\beta$ and $\scapro{\eta}{f} \ge \beta$ for all $\eta\in M$. Then
$\scapro{\theta}{f} =\frac{1}{t}\int_0^t \scapro{\zeta(s)}{f}\, ds \ge \beta$. Thus
$\theta\in M$ and therefore $\theta\in M\cap \mathcal{P}$.

Since $\{\zeta(s):s\ge 0\}$ is tight, its convex hull is tight and hence relatively
compact in $\calP$ by Prohorov's Theorem.
Thus the set $\{\theta_t:t\ge 0\}$ is contained in a compact set and therefore there
exist a sequence
 $t_n\uparrow \infty$ and a measure $\eta\in\calP$ such that
$\theta_{t_n}\to\eta$.

Finally, for $t\ge \alpha$ and $f\in C_b$ we have
\begin{align*}
\lim_{n\to\infty} \frac{1}{t_n}\int_0^{t_n} \langle \zeta(t+s),f\rangle ds
&=\lim_{n\to\infty} \frac{1}{t_n} \int_0^{t_n} \langle T_{s}^\ast,P_tf\rangle ds =
\langle \eta,P_tf\rangle
\intertext{and on the other hand}
\frac{1}{t_n}\int_0^{t_n} \langle \zeta(t+s),f\rangle ds
 &= \frac{1}{t_n} \int_0^{t_n} \langle \zeta(s),f\rangle ds
  -\frac{1}{t_n} \int_0^t \langle \zeta(s),f\rangle ds \\
&\qquad\qquad + \frac{1}{t_n}\int_{t_n}^{t_n+t} \langle \zeta(s),f\rangle ds,
\end{align*}
which converges to $\langle\eta,f\rangle$ as $n\to\infty$. Hence,
$P_t^\ast\eta=\eta$ and it follows that $P_t^\ast\eta=
P_t^\ast(P_\alpha^\ast\eta)=P_{t+\alpha}^\ast\eta=\eta$, for every
$t\ge 0$.
\end{proof}

We remark that the proof given above remains true in a more general setting. Indeed, we
need only to replace $D[-\alpha,0]$ by an arbitrary separable metric space $E$ and assume
that $(T_t)_{t\ge 0}$ is a family of bounded linear operators on $C_b(E)$ such that
$T_{s+t+\alpha}=T_sT_t$ for all $s,t\ge 0$,  that for some $\zeta\in\calP(E)$ one has
that $T_t^\ast\zeta\in\calP(E)$ for all $t\ge 0$, that the map $t\mapsto \langle
T_t^\ast\zeta,f\rangle$ is measurable from $[0,\infty)$ to $\R$ for all $f\in C_b(E)$,
and that the set $\{T_t^\ast\zeta:t\ge 0\}$ is tight.

\section{Uniqueness of the stationary solution}

As we have seen, the Markovian semigroup is in general not
eventually strongly Feller so that a main tool to establish
uniqueness of the invariant measure is not available. Moreover, when
considered as a stochastic evolution equation, the generator of the
deterministic equation \eqref{eq.det} is only eventually compact
(see \cite{DGVW95}) and the Markov semigroup is only weakly
continuous with a generator which is analytically not easily
tractable (see \cite{Mohammed84}). Hence, typical analytical methods
to prove uniqueness (see \cite{MaSei98} for a survey) cannot be
easily applied either.

We therefore consider several specific cases where uniqueness can be
proved nevertheless: for small Lipschitz constants by a contraction
argument, in the Wiener case for non-delayed diffusion coefficients
by establishing the strong Feller property via Girsanov's theorem
and for compound Poisson driving processes and non-delayed drift
terms by studying the deterministic behaviour between the jumps.
After that, we relax the requirements and show that in full
generality second-order uniqueness holds up to a constant factor. We
conclude by an example where the invariant measures are not unique.

\subsection{Small Lipschitz constants}

If the function $F$ is not too far from being constant, as
measured by the Lipschitz constant, then uniqueness holds.
The upper bound for the Lipschitz constant below can be reconstructed by
our proof, but it is certainly not the best possible.

\begin{theorem}
Grant Assumption \ref{AssStat} and suppose that the L\'evy process has finite second
moments. If the Lipschitz-constant $K$ of $F$ in \eqref{eq.lipcond} is sufficiently small
then the laws of all stationary solutions $X$ of \eqref{eq.stoch} coincide.
\end{theorem}

\begin{proof}
Let $X$ and $Y$ be two stationary solutions with corresponding initial
conditions  $X_0$
and $Y_0$. As mentioned above Proposition \ref{pro.second-order} the moments
$\E\norm{X_0}^2_\infty$ and $\E\norm{Y_0}^2_\infty$ are finite.

As $v_0(\mu)<0$ the fundamental solution $r$ decays exponentially with $\abs{r(t)}\le
ce^{-\beta t}$, $\int_0^\infty r^2(s)\exp(2\beta s)\,ds<\infty$, and $\int_0^\infty
\dot{r}^2(s)\exp(2\beta s)\,ds<\infty$ for some constants  $c,\,\beta>0$. We choose an
arbitrary constant $\gamma<\beta$ and we let $ Z(u):=F(X)(u)-F(Y)(u)$ for convenience. By
use of the decomposition $L(t)=M(t)+\E[L(1)]\,t$ with a martingale $M$, the variation of
constants formula implies for $t\ge \alpha$:
\begin{align}\label{eq.Ldecest}
\begin{split}
&\E\left[\sup_{t-\alpha\le s\le t}\abs{e^{\gamma s} (X(s)-Y(s))}^2\right]
 \le 3\E\left[\sup_{t-\alpha\le s\le t}\abs{e^{\gamma s}x(s,X_0-Y_0)}^2\right]\\
&\qquad + 3 \E\left[\sup_{t-\alpha\le s\le t}\abs{\int_0^s e^{\gamma s}
      r(s-u)Z(u-)\,dM(u)}^2\right] \\
&\qquad + 3(\E[L(1)])^2 \E\left[\sup_{t-\alpha\le s\le t}\abs{\int_0^s e^{\gamma
s} r(s-u)Z(u-)
\,du}^2\right] .
\end{split}
\end{align}
An application of representation \eqref{eq.xdet} yields
\begin{align}
 \E\left[\sup_{t-\alpha\le s\le t}\abs{e^{\gamma s}x(s,X_0-Y_0)}^2\right]
   \le d \E\norm{X_0-Y_0}_\infty^2         \label{eq.1est}
\end{align}
for a finite constant $d$ depending only on the measure $\mu$. Let $r_1$ be the
function
defined by $r_1(s):=r(s)\exp(\gamma s)$. Then we obtain for the second term in
\eqref{eq.Ldecest}
\begin{align}
&\E\left[\sup_{t-\alpha\le s\le t}\abs{\int_0^s e^{\gamma s}
      r(s-u)Z(u-)\,dM(u)}^2\right]\notag\\
&\qquad =\E\left[\sup_{t-\alpha\le s\le t}\abs{\int_0^s \left(r_1(0)+\int_0^{s-
u}\dot{r}_1(m)\,dm\right)
 e^{\gamma u}Z(u-)\,dM(u)}^2\right]\notag\\
 \begin{split}
&\qquad\le 2 \E\left[\sup_{t-\alpha\le s\le t}\abs{\int_0^s e^{\gamma u}Z(u-
)\,dM(u)}^2\right]\\
&\qquad\qquad  + 2\E\left[\sup_{t-\alpha\le s\le t}\abs{\int_0^s
\left(\int_0^{s-m} e^{\gamma u} Z(u-)\,dM(u)\right)
       \dot{r}_1(m)\,dm}^2\right]. \label{eq.rdotr}
\end{split}
\end{align}
The first term in \eqref{eq.rdotr} can be estimated by
\begin{align}
&\E\left[\sup_{t-\alpha\le s\le t}\abs{\int_0^s e^{\gamma u}Z(u-
)\,dM(u)}^2\right]\notag \\
&\qquad\qquad\le  \left(\sigma^2+\int x^2\,\nu(dx)\right)\int_0^t
  e^{2\gamma u}\E\abs{Z(u-)}^2\,du. \label{eq.2est}
\end{align}
Note that $\dot{r}_1$ has essentially the same asymptotic as $s\mapsto
\exp(\gamma s)r(s)$.
Hence if we choose a constant $\delta>0$ such that $\gamma + \delta \le \beta$
we obtain
$d_1:=\int_0^\infty \exp(2\delta m)\abs{\dot{r}_1(m)}^2\,dm<\infty$.
Applying H{\"o}lder's inequality to the second term in
\eqref{eq.rdotr} results in
\begin{align}
& \E\left[\sup_{t-\alpha\le s\le t}\abs{\int_0^s \left(\int_0^{s-m} e^{\gamma u}
Z(u-)\,dM(u)\right)
       \dot{r}_1(m)\,dm}^2\right]\notag\\
& \le \E\left[\sup_{t-\alpha\le s\le t}\int_0^s e^{2\delta
m}\abs{\dot{r}_1(m)}^2\,dm\;\;
 \!\! \int_0^s \abs{\int_0^{s-m} e^{\gamma u}Z(u-)\,dM(u)}^2\!\! e^{-2\delta
m}\,dm\right]    \notag \\
& \le d_1 e^{-2\delta (t-\alpha)} \E\left[\sup_{t-\alpha\le s\le t}\!
     \int_0^s \abs{\int_0^{m} e^{\gamma u}Z(u-)\,dM(u)}^2
      e^{2\delta m}\,dm\right]\notag\\
&\le d_1 e^{2\delta\alpha}\left(\sigma^2 +\int
x^2\,\nu(dx)\right)\left(\int_0^\infty e^{-2\delta m}\,dm\right)
     \int_0^t e^{2\gamma u}\E\abs{Z(u-)}^2\,du.\label{eq.3est}
\end{align}
The last term in \eqref{eq.Ldecest} can be estimated similarly by H{\"o}lder's
inequality
\begin{align}
&\E\left[\sup_{t-\alpha\le s\le t}\abs{\int_0^s e^{\gamma s} r(s-u)Z(u-)
\,du}^2\right] \notag\\
& \qquad \le \left(\int_0^\infty e^{2\gamma u}\abs{r(u)}^2 \,du\right)
 \int_0^t e^{2\gamma u} \E\abs{Z(u-)}^2\,du. \label{eq.4est}
\end{align}
By collecting the inequalities \eqref{eq.1est} to \eqref{eq.4est}, using
the Lipschitz condition \eqref{eq.lipcond} and applying Gronwall's Lemma we conclude
$\E\norm{X_s-Y_s}^2_\infty \to 0$ for $s\to \infty$ if the
Lipschitz constant $K$ is sufficiently small. Consequently, the laws of $X_0$ and $Y_0$ coincide.
\end{proof}

\subsection{Non-delayed diffusion coefficient}\label{sec.nondelaydiff}

We have seen in Section \ref{sse.feller} that the Markov semigroup
$(P_t)_{t\ge 0}$ of the solution segments is in general not
eventually strong Feller. This is only an effect due to the delay in
the diffusion term and cannot be caused by a delayed drift for the
Wiener-driven case, as we shall see now. Let us consider as special
case of equation \eqref{eq.stoch}
\begin{equation}\label{eq.nodelaydif}
dX(t)=\left(\int_{[-\alpha,0]}X(t+s)\,\mu(ds)\right)\,dt+f(X(t))\,dW(t)\quad
\text{ for
}t\ge 0,
\end{equation}
with initial segment $\Phi$ as in \eqref{eq.stoch}, a Wiener process $W$ and a
Lipschitz
function $f:\R\to \R$. By a simple argument based on Girsanov's theorem we
obtain the
following result.

\begin{proposition}
If $f$ satisfies the ellipticity condition $\inf_{x\in\R}f(x)>0$, then the
solution
segments $(X_t:\, t\ge 0)$ of \eqref{eq.nodelaydif} generate a Markov semigroup
on
$C([-\alpha,0])$ that is strongly Feller after time $\alpha$.
\end{proposition}

\begin{proof}
First note that the continuous functions form a closed subspace of
the Skorokhod space $D[-\alpha,0]$ such that the formerly obtained
results are in the Wiener-driven case also valid on $C[-\alpha,0]$.
Referring to Theorem 7.19 for diffusion-type processes in
\cite{LipShi01}, we infer from the Lipschitz continuity of the
coefficients and from the ellipticity of $f$ that the laws $Q_1$ and
$Q_2$ of the solution processes of \eqref{eq.nodelaydif} on
$C[0,T]$, $T>0$ arbitrary, are equivalent for different delay
measures $\mu_1$ and $\mu_2$. The corresponding Radon-Nikodym
derivative is given by
\begin{align*}
\frac{dQ_2}{dQ_1}(X)&=\exp\biggl(\int_0^T
\left(\int_{[-\alpha,0]}X(t+s)\,(\mu_1-\mu_2)(ds)\right)
f(X(t))^{-2}\,dX(t)\\
&\qquad-\frac{1}{2}\int_0^T \left(\int_{[-\alpha,0]}X(t+s)\,(\mu_1-
\mu_2)(ds)\right)^2
f(X(t))^{-2}\,dt\biggr).
\end{align*}
As in \cite[Thm. 2.1]{MasSei00} one can show that the validity of the strong
Feller
property at each time is invariant under the change of measure. According to
that result
we need to check that the semigroup is Feller and that
\[
\lim_{n\to\infty}\E\abs{\frac{dQ_2}{dQ_1}(X^1(\cdot;\phi^n))
-\frac{dQ_2}{dQ_1}(X^1(\cdot;\phi))}=0
\]
for initial segments $\phi^n\to\phi$ in $C[-\alpha,0]$ and for the
corresponding solution process $X^1$ with the choice $\mu_1$. The
Feller property has been established in Section \ref{sse.feller}.
By Scheff\'e's Lemma it suffices for the second condition to prove
convergence in probability. This is accomplished by the continuity
of the map $\phi\mapsto X^1(\cdot,\phi)$ from $C[-\alpha,0]$ to
$L^2([0,T]\times\Omega)$ for any $T$, which follows from
\cite[Thm. 3.1]{Mohammed84}.

We have thus reduced the problem to proving the strong Feller property of the
Markov
semigroup generated by the solution segments $(\tilde{X}_t)_{t\ge 0}$ of
\begin{equation}\label{eq.nodelaydif2}
d\tilde{X}(t)=f(\tilde{X}(t))\,dW(t)\quad \text{ for }t\ge 0,
\end{equation}
as special case of \eqref{eq.nodelaydif} with $\mu=0$. It is well known that
this
diffusion equation generates a strongly Feller semigroup on $\R$ under our
assumptions on $f$,
see e.g. \cite[Thm.7.1.1]{DaPZab96}. We claim that this property is
inherited
by the segment process. For this consider a bounded measurable functional $\Psi$
on
$C[-\alpha,0]$ and remark that $\tilde{X}(\cdot;\phi)=\tilde{X}(\cdot;\phi(0))$
only
depends on the initial value, not the whole segment. By the scalar Markov and
weak
uniqueness property we obtain for $t\ge \alpha$ and any initial segment $\phi$
with
obvious notation
\[ \E[\Psi(\tilde{X}_t(\phi))]
=\E[\E[\Psi(\tilde{X}_t(\phi))\,|\,{\cal F}_{t-\alpha}]]
=\E_\omega[\E_{\omega'}[\Psi(\tilde{X}_\alpha(\tilde{X}(t-
\alpha;\phi,\omega),\omega'))]].
\]
Setting $H(\xi):=\E[\Psi(\tilde{X}_\alpha(\xi))]$, $\xi\in\R$, the scalar strong
Feller
property implies the continuity of
\[ \eta\mapsto \E[H(\tilde{X}(t-\alpha;\eta))]
 =\E_\omega[\E_{\omega'}[\Psi(\tilde{X}_\alpha(\tilde{X}(t-
\alpha;\eta,\omega),\omega'))]]
\]
for $\eta\in\R$. Since $\phi^n\to\phi$ in $C[-\alpha,0]$ yields $\phi^n(0)\to
\phi(0)$,
we thus infer the continuity of
\[ \phi\mapsto \E_\omega[\E_{\omega'}[\Psi(\tilde{X}_\alpha(\tilde{X}(t-
\alpha;\phi,\omega)),\omega')]]
=\E[\Psi(\tilde{X}_t(\phi))]
\]
on $C[-\alpha,0]$, which is the asserted strong Feller property at $t\ge
\alpha$.
\end{proof}

\begin{corollary}
The Markov semigroup $(P_t)_{t\ge 0}$ is regular after time
$2\alpha$. Thus, any stationary solution of \eqref{eq.nodelaydif} is unique and
strongly
mixing.
\end{corollary}

\begin{proof}
Recall that we have regularity at $t_0$ if all transition probabilities
$P(X_{t_0}(\phi)\in\cdot)$ are equivalent for $\phi\in C[-\alpha,0]$. By Doob's
Theorem
\cite[Thm. 4.2.1]{DaPZab96} this property yields the uniqueness and strong
mixing result.

The regularity property at $t_0>2\alpha$ is implied by the strong Feller
property at
$\alpha$ together with the irreducibility at $t_0-\alpha$ \cite[Prop.
4.1.1]{DaPZab96},
which means that all transition probabilities at time $t_0-\alpha$ have support
in the
entire space. To prove the latter, we may again restrict to the case $\mu=0$ and
consider
$\tilde{X}$ as in \eqref{eq.nodelaydif2} due to the equivalence of the laws. As
in
\cite[Cor. VIII.2.3]{RevYor99} it follows from Girsanov's theorem that the
support of the
(regular) conditional law $\law(\tilde{X}_{t_0-\alpha}\,|\,\tilde{X}(t_0-
2\alpha)=x)$ is
given by
$ S_x:=\{f\in C[-\alpha,0]\,:\,f(-\alpha)=x\}.$
Since the law of $\tilde{X}(t_0-2\alpha;\phi)$ has for the same reasons the full
support
$\R$ for any initial segment $\phi$, we conclude by composition that
$\law(X_{t_0-\alpha}(\phi))$ has full support $C[-\alpha,0]$ independent of
$\phi$,
which yields the required irreducibility.
\end{proof}

\subsection{Uniqueness in the compound Poisson case}

Let us consider here the case of a L\'evy triplet $(b,\sigma^2,\nu)$
with $\sigma=0$, $b=0$, and the total variation
$\lambda:=\norm{\nu}_{TV}$ finite, that is $L$ is a compound Poisson
process. If there is no delay in the drift, then we can reduce the
question of uniqueness of the invariant law on the Skorohod space
$D[-\alpha,0]$ to a property of the one-dimensional invariant law.

\begin{proposition}\label{prop.CPuniq}
Suppose $L$ is a compound Poisson process and consider for $a>0$ a
differential equation of the form
\begin{equation}\label{eq.adL}
dX(t)=-aX(t)\,dt + F(X)(t-)\,dL(t)\quad\text{for } t\ge 0,
\end{equation}
admitting a strong solution for any initial segment. If an invariant solution
measure on
$D[-\alpha,0]$ exists and the one-dimensional marginal distributions of any two
invariant
measures are non-singular, then the invariant measure is unique.
\end{proposition}

\begin{proof}
Let $\rho_1$ and $\rho_2$ be two invariant measures. By coupling methods we can construct
a filtered probability space $(\Omega,{\cal F},({\cal F}_t),P)$ carrying the process $L$,
and the ${\cal F}_0$-measurable random variables $Y,Z\in D[-\alpha,0]$ with $P^Y=\rho_1$,
$P^Z=\rho_2$ and $P(Y(0)=Z(0))>0$. Denote by $X^1$ and $X^2$ the corresponding strong
solution processes with initial conditions $Y$ and $Z$, respectively.

Since with probability $e^{-\lambda\alpha}>0$ the process $L$ does
not jump on the interval $[0,\alpha]$, we have
\begin{align*}
P(X^1_\alpha=X^2_\alpha,\, Y(0)=Z(0))
&\ge P(\sum_{t\le \alpha}\abs{\Delta L(t)}=0, Y(0)=Z(0))\\
& =P(\sum_{t\le \alpha}\abs{\Delta
L(t)}=0)P(Y(0)=Z(0))>0.
\end{align*}
Hence, introducing the set
\[S:=\{\phi\,|\,\exists\omega\in\Omega:\:X^1_\alpha(\omega)=X^2_\alpha(\omega)=\phi\}
 \subseteq D[-\alpha,0],
\]
we find for any Borel set $B$ in $D[-\alpha,0]$
\begin{align*}
P(X^1_\alpha\in B\cap S)&\ge
P(\{\omega\,|\,X_\alpha^2(\omega)=X_\alpha^1(\omega),\,
X_\alpha^1(\omega)\in
B\})\\
&=P(\{\omega\,|\,X_\alpha^2(\omega)=X_\alpha^1(\omega),\,
X_\alpha^2(\omega)\in B\})\\
&=:\tau_S(B)
\end{align*}
and equivalently $P(X^2_\alpha\in B\cap S)\ge \tau_S(B)$. By
invariance, we conclude
\[ \min\{\rho_1(B\cap S),\rho_2(B\cap S)\}\ge\tau_S(B)\quad\text{for all } B\in
{\cal B}(D[-\alpha,0])
\]
with a non-negative measure $\tau_S$ satisfying $\tau_S(S)>0$. Hence $\rho_1$ and
$\rho_2$ are non-singular; for if $\rho_1(A)=0$ and $\rho_2(A^C)=0$ for some Borel set
$A$, then
\[ \tau_S(S)=\tau_S(S\cap A)+\tau_S(S\cap A^C)\le \rho_1(S\cap
A)+\rho_2(S\cap A^C)=0.
\]
As extremal points of the set of invariant measures are singular (see \cite[Prop.
3.2.7]{DaPZab96}), uniqueness follows.
\end{proof}

\begin{theorem}
Grant Assumption \ref{AssStat}. Suppose $L$ is a compound Poisson process and
consider
equation \eqref{eq.adL} with $a>0$ and $F(\phi)(0)>0$ for all $\phi\in
D[-\alpha,\infty)$. Then there exists a unique invariant measure for
\eqref{eq.adL}.
\end{theorem}

\begin{proof}
If the jump measure $\nu$ is zero, then the compound Poisson
process vanishes and the only invariant measure is clearly the
point measure in zero. Let us now first consider the case of
possible positive jumps: $\nu((0,\infty))>0$. By Proposition
\ref{prop.CPuniq} it suffices to show that any two invariant
one-dimensional distributions are non-singular. We first show that
they are absolutely continuous with respect to Lebesgue measure.

Let $B\subseteq\R$ denote any Borel set. For the solution process $X$ of
\eqref{eq.adL}
  we find
\begin{align}
P(X(t)\in e^{-at}B)&\ge P\Bigl(\sum_{s\le t}\abs{\Delta
L(s)}=0,\, X(0)e^{-at}\in e^{-at}B\Bigr)\nonumber\\
&=e^{-\lambda t}P(X(0)\in B).\label{eq.mu0est}
\end{align}
Now assuming that $X$ is stationary with one-dimensional marginal
law $\rho_0$, we obtain by Fubini's Theorem for any Lebesgue null
set $B$ and $T>0$
\begin{align*}
\int_0^T\rho_0(e^{-at}B)\,dt&=\int_0^T\int_{\R}
\1_{e^{-at}B}(x)\,\rho_0(dx)\,dt\\
&=\int_{\R\setminus\{0\}}\int_x^{xe^{aT}}\frac{a}{t} \1_{B}(t)\,dt
       \,\rho_0(dx)+\rho_0(\{0\})\1_B(0)\\
&=\rho_0(\{0\})\1_B(0).
\end{align*}
By estimate \eqref{eq.mu0est}, however, the left-hand side is bounded from
below by $\rho_0(B)\frac{1-e^{-\lambda T}}{\lambda}$. Hence, we infer
$\rho_0(B)=0$ for all Lebesgue null sets $B$ with $0\notin B$. Since $F$ is
positive and $\nu\not=0$, we can exclude a point mass in zero because the state
$\{0\}$ will be eventually left by the process $P$-a.s. and the probability to
jump back exactly to this state is zero. We conclude that $\rho_0$ is
absolutely continuous with respect to the Lebesgue measure.

Let $S$ denote the support of $\rho_0$. Since $F$ is positive and bounded away
from zero and $L$ has positive jumps, there will occur with positive
probability sufficiently many positive jumps of $L$ in short time that the
trajectory $X$ will take arbitrarily high values. This means for the support
$S$ of the marginal invariant measure $\rho_0$ that $\sup S=+\infty$.

For a Borel set $B\subseteq (0,\infty)$ we have
\begin{align*}
\int_0^\infty e^{-at}\rho_0(e^{at}B)\,dt
&=\int_0^\infty \int_{(0,\infty)} e^{-at}\1_{e^{at}B}(x)\,\rho_0(dx)\,dt\\
&= \int_{(0,\infty)} \int_0^x \frac{1}{ax}\1_{B}(s)\,ds\,\rho_0(dx).
\end{align*}
If $\rho_0(B)=0$, then (\ref{eq.mu0est}) with $B$ replaced by $e^{at}B$ yields
that $\rho_0(e^{at}B)\le e^{\lambda t}\rho_0(B)=0$ for all $t>0$, and we obtain
that
\[
\int_{(0,x)} \1_B(s)\,ds =0\quad\mbox{for }\rho_0\mbox{-a.e.}\ x>0.
\]
Since $\sup S=+\infty$, we infer that the Lebesgue measure of $B$ equals $0$.
Thus the Lebesgue measure of $(0,\infty)$ is absolutely continuous with respect
to $\rho_0$.

If $\nu((-\infty,0))$ is also positive, then the symmetric argument
yields that the Lebesgue measure on $\R$ is equivalent with
$\rho_0$. In any case, we know that two invariant measures are both
equivalent to the appropriate Lebesgue measure and hence with each
other. An application of Proposition \ref{prop.CPuniq} completes the
proof.
\end{proof}

\begin{remark}

We have derived the regularity property that, unless the
jump measure is zero, the one-dimensional marginals are absolutely
continuous with respect to the Lebesgue measure.

In some cases one can easily derive the density of the invariant measure. For example, if
we assume $L$ to have only positive jumps of size at least $J>0$ and $F(\phi)(0)\in
[\sigma_0,\sigma_1]$ for all $\phi\in D[-\alpha,\infty)$ and some $\sigma_0,\,\sigma_1>0$
then the density of the marginal of the invariant measure of \eqref{eq.adL} near zero is
given by
\begin{align*}
f(x)=C\frac{\lambda}{\alpha}x^{(\lambda-a)/a},
\quad x\in [0,J\sigma_0),
\end{align*}
with a suitable constant $C$.

\end{remark}

\subsection{Second-order uniqueness}
A real-valued stochastic process $(X(t):t\ge -\alpha)$ will be
called {\em second-order stationary}, if $0<\E[X(t)]^2<\infty$, the
values $\E[X(t)]$ are constant for all $t\ge -\alpha$, and the
function $(s,t)\mapsto \E[X(s)X(t)]$
depends only on the difference $s-t$. Obviously, any stationary
 solution of \eqref{eq.stoch} with finite second moments is
second-order stationary. If the L\'evy process is a
square-integrable martingale, we establish second-order uniqueness
for equation \eqref{eq.stoch} up to a constant factor, more
precisely the expectation and the correlation function are
uniquely determined and can be calculated analytically.

Note that the invariant measure exhibited in Section
\ref{sec.statexist} will have finite second moments for its
one-dimensional marginal whenever the L\'evy process has finite
second moments. This follows from the fact that the constructed
tight sequence of segments $(X_t)$ will be uniformly bounded in
$L^2_P(\Omega)$ by the variation of constants formula
\eqref{eq.varcons} and Lemma \ref{lem.emery}.

\begin{proposition}\label{pro.second-order}
Grant Assumption \ref{AssStat}. Suppose the L\'evy process is a
square-integrable martingale with characteristics
$(b,\sigma^2,\nu)$. Then any stationary solution $(X(t):\,t\ge
-\alpha)$  of \eqref{eq.stoch} with finite second moments is a
centered random process with auto-covariance function
\begin{align*}
c(h):=\E[X(0)X(h)]=\frac{\Var[X(0)]}{\norm{r}_{L^2(\Rp)}^2}\int_0^\infty
r(s)r(s+h)\,ds,\quad h\ge 0.
\end{align*}
The spectral density is given by
\begin{align*}
\xi\mapsto
\E[X(0)^2]\left(\norm{r}_{L^2(\Rp)}\abs{\chi_\mu(i\xi)}\right)^{-2},
\qquad \xi\in\R,
\end{align*}
 where $\chi_\mu(z):=z-\int_{[-\alpha,0]}e^{zu}\,\mu(du)$ is the
characteristic function of the deterministic equation
\eqref{eq.det}.
\end{proposition}

\begin{proof}
By the variation of constants formula \eqref{eq.varcons} and the
martingale property of $L$ we have for $t\ge 0$
\[ \E X(t)=\E  x(t,X_0)=r(t)\E[X(0)]+\int_{[-\alpha,0]}\int_{s}^0
r(t+s-u)\E[X(u)]\,du\,\mu(ds).
\]
Due to $\lim_{t\to\infty} r(t)=0$ and stationarity we conclude
that $X$ is centered. Again using the variation of constants
formula, we find for $h,\,t\ge 0$
\begin{align*}
&\E[X(t)X(t+h)] =  \E[x(t+h,X_0)\int_0^t r(t-u)F(X)(u-)\,dL(u)]\\
& +\E [x(t,X_0)\int_0^{t+h} r(t+h-u)F(X)(u-)\,dL(u) ]
    + \E\left[x(t,X_0)\;x(t+h,X_0)\right]\\
& + \E[\int_0^t r(t-u)F(X)(u-)\,dL(u)
            \int_0^{t+h} r(t+h-u)F(X)(u-)\,dL(u)].
\end{align*}
As in Lemma \ref{lem.emery} we obtain
\begin{align*}
& \E\left[\int_0^t r(t-u)F(X)(u-)\,dL(u)\int_0^{t+h}
 r(t+h-u)F(X)(u-)\,dL(u)\right]\\
&\qquad= \left(\sigma^2 + \int x^2 \nu(dx)\right)\int_0^t
 r(t-u)r(t+h-u)\E[F(X)(u-)]^2\,du.
\end{align*}
The variance is estimated as the expectation before:
\begin{align*}
&\Var[x(t,X_0)]\\
&\;\; \le 2
\left(r(t)^2\Var[X(0)]+\left(\int_{[-\alpha,0]}\int_{-s}^0
\abs{r(t+s-u)}\E|X(u)|\,du\,\abs{\mu}\!(ds)\right)^2\right),
\end{align*}
which converges to $0$ as $t\to\infty$. Applying the
Cauchy-Schwarz inequality to the first three terms in the equation
above results in
\begin{align*}
 \Cov(X(0),X(h)) &=\lim_{t\to\infty}\E[X(t)X(t+h)]\\
 &= \E[F(X)(0)]^2 \left(\sigma^2 + \int x^2 \nu(dx)\right)
  \int_0^\infty r(u)r(u+h)\,du.
\end{align*}
This yields the expression for the covariance function. The
formula for the spectral density follows from the fact that $r$ is
the inverse Fourier transform of $\chi_\mu(-i\cdot)^{-1}$, as
obtained for affine stochastic delay differential equations driven
by a Wiener process in \cite{KuMe92}.
\end{proof}

\begin{remark}\label{re.second-order}
It is seen from the proof that
\[ \Var[X(0)]=\E[ F(X)(0)]^2 \left(\sigma^2 + \int x^2 \nu(dx)\right)
\norm{r}_{L^2(\Rp)}^2,
\]
which gives some information about the size of the variance
depending on bounds for the functional $F$. We shall see in the
counterexample of Section \ref{sec.nonuniq} that this variance
term need not be uniquely determined, at least for measurable
functionals $F$.
\end{remark}

\subsection{Non-uniqueness}\label{sec.nonuniq}

In the Wiener-driven case we construct an elliptic diffusion
functional $F$ which remains constant in time for certain initial
segments, but with different values for different initial segments.
By doing so, we can recover, for instance, Ornstein-Uhlenbeck
processes with different diffusion coefficients as solutions.
Suppose $F$ is of the form
\begin{align*}
 F(\phi)(t):= \sqrt{\max(1,\, \min(\tfrac{2}{\alpha}\langle \phi\rangle_{t-
\alpha}^{t-\alpha/2},\,2))}
  \1_{\Rp}(t)
\quad\text{for }t\ge -\alpha,
\end{align*}
where $\langle \phi\rangle_a^b$ denotes the quadratic variation of
$\phi\in D[-\alpha,\infty)$ on the interval $[a,b]$ which might be
infinite. Then $F$ is bounded away from zero and infinity and is
measurable (as a limit of measurable functionals), but obviously
not continuous. Leaving slightly our framework, let us consider
for a Wiener process $W$ the equation
\begin{align}
\begin{split}
 dX(t)&= -X(t)\,dt + F(X)(t-
)\,dW(t)\quad\text{for }t\ge 0,\\
  X(u)&= \Phi(u)\quad\text{for } u\in [-\alpha,0].
 \end{split}\label{SDDECE}
\end{align}
Due to the positive minimal delay $\alpha/2$ there exists a strong
unique solution to this equation for any ${\cal F}_0$-measurable
initial segment by the method of steps, cf. Mao \cite{Mao97}. On
the other hand, there exists for every $\sigma \in [1,2]$ a
stationary Ornstein-Uhlenbeck process $X^\sigma$ which solves the
equation (we suppose that $W$ is a two-sided Wiener process)
\begin{align*}
dX^\sigma(t)=-X^\sigma(t)\,dt + \sigma\,dW(t) \quad \text{for
}t\in\R .
\end{align*}
Then choosing $\Phi^\sigma=X^\sigma_0$, we obtain that each
$X^\sigma$ is also a stationary solution of \eqref{SDDECE}. This
is due to the fact that $\langle
X^\sigma\rangle_{t-\alpha}^{t-\alpha/2}=\frac{\alpha}{2}\sigma^2$
and thus $F(X^\sigma)(t)=\sigma$ hold for all $t\ge 0$ and
$\sigma\in [1,2]$.

This example shows that some kind of regularity of $F$ has to be
imposed to guarantee uniqueness, but we do not know whether
already for functionals $F$ with large, but finite Lipschitz
constants uniqueness breaks down. It is interesting to note that a
similar dichotomy has been described by Mohammed and Scheutzow
\cite{MohSch97} for the long time behaviour in dependence of the
diffusion functional.

\bibliographystyle{plain}

\bibliography{litstationary}

\end{document}